\documentclass[12pt]{article}
\usepackage{amsmath,amssymb,amsthm,enumerate,nicefrac,color,fancyhdr,hyperref,graphicx,adjustbox}

\usepackage[T1]{fontenc}

\usepackage{algorithm,caption}
\usepackage{algpseudocode}
\usepackage[citestyle=numeric,bibstyle=numeric,maxbibnames=99]{biblatex}
\usepackage{tkz-graph}
\usetikzlibrary{decorations.markings}
\tikzset{middlearrow/.style={
        decoration={markings,
            mark= at position 0.5 with {\arrow{#1}} ,
        },
        postaction={decorate}
    }
}
\renewcommand{\proof}[1]{\noindent\textbf{Proof#1.}}

\newcommand{\ceil}[1]{\left\lceil #1 \right\rceil}
\renewcommand{\qed}{\hfill\blacksquare}
\usepackage[margin=3cm]{geometry}
\newtheorem{theorem}{Theorem}[section]

\newtheorem{stheorem}{Theorem}[subsection]
\newtheorem{slemma}[stheorem]{Lemma}
\newtheorem{sdef}[stheorem]{Definition}
\newtheorem{problem}{Problem}
\newtheorem{prop}{Proposition}

\newtheorem{scorollary}[stheorem]{Corollary}
\newtheorem{lemma}[theorem]{Lemma}
\newtheorem{innercustomgeneric}{\customgenericname}
\providecommand{\customgenericname}{}
\newcommand{\newcustomtheorem}[2]{%
  \newenvironment{#1}[1]
  {%
   \renewcommand\customgenericname{#2}%
   \renewcommand\theinnercustomgeneric{##1}%
   \innercustomgeneric
  }
  {\endinnercustomgeneric}
}
\usetikzlibrary{decorations.pathreplacing,angles,quotes}
\newcustomtheorem{customthm}{Theorem}
\newcustomtheorem{customlemma}{Lemma}
\newcustomtheorem{customcor}{Corollary}
\title{Complexity of the Feedback Vertex Set Problem in Tournaments with Forbidden Subtournaments}
\author{Sophie Spirkl\thanks{We acknowledge the support of the Natural Sciences and Engineering Research Council of Canada (NSERC), [funding reference number RGPIN-2020-03912].
    Cette recherche a \'et\'e financ\'ee par le Conseil de recherches en sciences naturelles et en g\'enie du Canada (CRSNG), [num\'ero de r\'ef\'erence RGPIN-2020-03912]. This project was funded in part by the Government of Ontario. This research was conducted while Spirkl was an Alfred P. Sloan Fellow.}, Yun Xing\\Department of Combinatorics and Optimization, University of Waterloo}
\date{\today}
\begin{document}
\maketitle
\begin{abstract}
In this paper, we consider the complexity of the minimum feedback vertex set problem (MFBVS) for tournaments with forbidden subtournaments. The MFBVS problem in general tournaments is known to be NP-complete. We prove that the MFBVS problem for $W_5$-free and  $U_5$-free tournaments is  in P, and for $T_5$-free tournaments it remains NP-complete. Moreover, we prove a necessary condition for all $H$ such that the MFBVS problem for $H$-free tournaments is in P. We also show that the necessary condition is not sufficient.
\end{abstract}
\section{Introduction}
\subsection{Preliminaries}
A \textit{tournament} is a simple directed graph where for any pair of distinct vertices $(u,v)$, there is exactly one arc with ends $\{u,v\}$. If $u,v$ are vertices in a tournament and the arc between them starts at $u$ and ends at $v$, we denote this arc as $(u,v)$. In this paper, all tournaments are finite.\\

Here are some definitions used in this paper. Let $G$ be a tournament. A \textit{subtournament} of $G$ is an induced subgraph of $G$. Let $u,v \in V(G)$ and $X,Y \subseteq V(G)$. We use $u \rightarrow v$ to denote that $v$ is an out-neighbour of $u$, and $X \Rightarrow Y$ if each vertex in $X$ is an out-neighbour of all vertices in $Y$. We use $N^-(v)$ to denote the set of in-neighbors of $v$ and $N^+(v)$ to denote the set of out-neighbors of $v$. A tournament is \textit{transitive} if for every $u,v,w \in V(G)$ such that $u \rightarrow v$, $v \rightarrow w$, we have $u \rightarrow w$. A \textit{cyclic triangle} is the 3-vertex tournament that is not transitive. Motivated by the notion of coloring in directed graphs introduced by Neumann-Lara\cite{NEUMANNLARA1982265}, we call a set $X \subseteq V(G)$ an \textit{independent set} of $G$ if the induced subgraph $G[X]$ is transitive. We define the complement $\bar{G}$ of a tournament $G$ as the directed graph obtained by reversing every arc of $G$. We also define the complement of a vertex set $S \subseteq V(G)$ as $V(G) \setminus S$.\\

Based on \cite{BERGER20131}, we define the backedge graph of a tournament as follows. Let \( G \) be a tournament on \( n \) vertices, and let \( (v_1, \ldots, v_n) \) be an ordering of the vertex set \( V(G) \). An edge \( (v_i, v_j) \in E(G) \) is called a \emph{backedge} if \( i > j \). The \emph{backedge graph} \( B \) of $G$ is the undirected graph with vertex set \( \{v_1, \ldots, v_n\} \), where an edge \( v_i v_j \in E(B) \) exists if and only if \( (v_i, v_j) \) is a backedge in \( G \). We have the following well-known property for all tournaments:
\begin{prop}[Folklore]\label{folklore}
Let $G$ be a tournament. The following are equivalent:\\
(1) $G$ is transitive.\\
(2) $G$ contains no cyclic triangle.\\
(3) $G$ is acylic.\\
(4) $G$ admits a topological ordering on its vertices $V(G)$.\\
(5) There exists an ordering of $V(G)$ that gives a backedge graph with no edges.
\end{prop}

\begin{sdef}[Feedback vertex set]
Let $G$ be a tournament. A set $S \subseteq V(G)$ is a feedback vertex set if $V(G) \setminus S$ is an independent set.
\end{sdef}
\noindent In this paper, we are interested in finding the minimum feedback vertex set of a tournament. This problem is introduced by Speckenmeyer in 1990 and is NP-complete for general tournaments\cite{SPECKENMEYER}.
\begin{problem}[\textbf{MFBVS}: Minimum Feedback Vertex Set Problem for tournaments]
\ \\
Given a tournament $G$, find the minimum size of a feedback vertex set in $G$.
\end{problem}
\noindent Since the complement of a feedback vertex set is an independent set, the MFBVS problem is equivalent to the following problem:
\begin{problem}[\textbf{MISP}: Maximum Independent Set Problem for tournaments]
\ \\
Given a tournament $G$, find the maximum size of an independent set in $G$.
\end{problem}
\begin{problem}[\textbf{MISP, Decision Version}]
\ \\
Given a tournament $G$ and $k \in \mathbb{Z}_+$, does there exist an independent set in $G$ with size at least $k$?
\end{problem}

\noindent We can also define the weighted version of MFBVS as follows:
Let $G$ be a tournament and let $w: V(G) \rightarrow  \mathbb{Q}_{+}$ be a weight-assignment on the vertex set of $G$. The weight of a set $S \subseteq V(G)$ is defined as $\sum_{v \in S}w(v)$. Now we can define the weighted version of MFBVS:
\begin{problem}[\textbf{WMFBVS}: Weighted Minimum Feedback Vertex Set Problem for tournaments]
\ \\
Given a tournament $G$, find the maximum weight of a feedback vertex set in $G$.
\end{problem}
\noindent Again, this is equivalent to the following problem:
\begin{problem}[\textbf{WMISP}: Weighted Maximum Independent Set Problem for tournaments]
\ \\
Given a tournament $G$, find the maximum weight of an independent set in $G$.
\end{problem}
\begin{problem}[\textbf{WMISP, Decision Version}]
\ \\
Given a tournament $G$ and $k \in \mathbb{Q}_+$, does there exist an independent set in $G$ with weight at least $k$?
\end{problem}

\noindent In this paper, we mainly focus on studying MISP and WMISP instead of MFBVS and WMFBVS. We use MISP to refer Problem 3, the decision version, even though most of our algorithms will actually find a maximum independent set. Similarly, we use WMISP to refer Problem 5. Note that MISP is a special case of WMISP. 

\subsection{Motivation}
\noindent It is natural to ask whether MISP is NP-complete or not. The following theorem is due to Speckenmeyer:

\begin{theorem}[Speckenmeyer, \cite{SPECKENMEYER}]\label{npc}
MISP for tournaments is NP-complete.
\end{theorem}

\noindent We include a formal proof of this result, as \cite{SPECKENMEYER} only includes sketch of this proof. Our proof is a simplified version of Speckenmeyer's proof.\\
\proof{}
We first define the Minimum Vertex Cover Problem for graphs.\\
A \textit{vertex cover} $C \subseteq V(G)$ is a vertex set such that every $e \in E(G)$ has at least one end in $C$. A minimum vertex cover is a vertex cover of minimum size. 
\begin{problem}[\textbf{MVCP}]
    Given a graph $G$, find the size of a minimum vertex cover of $G$.
\end{problem}
It is well-known that MVCP is NP-complete. (A proof is given in Chapter 3 of \cite{Garey-Johnson}.) To show NP-completeness of MISP, we show that MVCP is reducible to MISP. Given a graph $G = (V,E)$, where $V = \{v_1,...,v_n\}$, we construct the following graph $G' = (V', E)$ where
\begin{equation*}
    V' = \{v_{1},v_{1}',...,v_{n},v_{n}'\}.\\
\end{equation*}

Let $T$ be a tournament on $V'$ whose backedge graph under the ordering $(v_{1},v_{1}',\allowbreak ...,v_{n},v_{n}')$ is $G'$. Let $C^*$ be a minimum vertex cover of $G$. We claim that $V' \setminus C^*$ is a maximum independent set in $T$.  Clearly, since $C^*$ is a minimum vertex cover in $G$, we have that $C^*$ is also a minimum vertex cover in $G'$. As a result, $G' \backslash C^*$ contains no edge, and then by Proposition \ref{folklore}, it follows that $T\backslash C^*$ is transitive, so $V' \setminus C^*$ is independent in $T$. Suppose for contradiction that there exists $S'$ such that $|S'| < |C^*|$ and $T \setminus S'$ is transitive. Now let $S = \{v_i \in V: \{v_{i},v_{i}'\} \cap S' \neq \varnothing, i \in [n]\}$. Suppose that $S$ is not a vertex cover of $G$, then there is an edge $e=\{v_i,v_j\} \in E(G\backslash S)$ with $i < j$. Now $\{v_{i},v_{i}',v_{j}\} \subseteq V' \setminus S'$ induces a cyclic triangle in  $T$, so by Proposition \ref{folklore}, $T\setminus S'$ is not transitive, which is a contradiction. Thus $S$ is a vertex cover of $G$. Since $|C^*| > |S'| \geq |S|$, this contradicts the minimality of $C^*$. Hence, finding the size of maximum independent set in $T$ gives the size of minimum vertex cover in $G$. $\qed$\\

\noindent Let $G,H$ be tournaments. We call $G$ \textit{H-free} if $G$ has no induced subgraph that is isomorphic to $H$. The main focus of this paper is to initiate the study of the following question: for which $H$ is WMISP restricted to \text{$H$-free} tournaments solvable in polynomial time (in terms of size of the tournament), and for what $H$ is MISP restricted to \text{$H$-free} tournaments NP-complete? This is motivated by the analogous question in graphs, see for example \cite{MUNARO, GIACOMO}.\\

Erd\H{o}s and Moser showed in \cite{erdos1964} that if $H$ is a transitive tournament, then the number of $H$-free tournaments is finite. Therefore WMISP can be solved in polynomial time for $H$-free tournaments if $H$ is transitive. If $H$ consists of 3 vertices and is not transitive (namely, $H$ is a cyclic triangle), then by Proposition \ref{folklore}, every $H$-free tournament is transitive. Thus, the first non-trivial case for our problem is when $|H| = 4$. Here are the four tournaments (up to isomorphism) with 4 vertices:\\

\begin{center}
\begin{tikzpicture}[thick,scale=0.6]
\node[fill,circle,inner sep=0pt,minimum size=3.5pt] (1) at (0,0){}; 
\node[fill,circle,inner sep=0pt,minimum size=3.5pt] (2) at (3,0){}; 
\node[fill,circle,inner sep=0pt,minimum size=3.5pt] (3) at (0,-3){}; 
\node[fill,circle,inner sep=0pt,minimum size=3.5pt] (4) at (3,-3){}; 
\draw[-stealth,semithick] (1) edge (2);
\draw[-stealth,semithick] (1) edge (3);
\draw[-stealth,semithick] (1) edge (4);
\draw[-stealth,semithick] (2) edge (3);
\draw[-stealth,semithick] (2) edge (4);
\draw[-stealth,semithick] (3) edge (4);
\node[circle, inner sep=1pt, label={$K_4$}] (a) at (1.5,-4.5) {};
\end{tikzpicture}\quad\quad
 \begin{tikzpicture}[thick,scale=0.6]
\node[fill,circle,inner sep=0pt,minimum size=3.5pt] (1) at (0,0){}; 
\node[fill,circle,inner sep=0pt,minimum size=3.5pt] (2) at (3,0){}; 
\node[fill,circle,inner sep=0pt,minimum size=3.5pt] (3) at (0,-3){}; 
\node[fill,circle,inner sep=0pt,minimum size=3.5pt] (4) at (3,-3){}; 
\draw[-stealth,semithick] (1) edge (2);
\draw[-stealth,semithick] (3) edge (1);
\draw[-stealth,semithick] (4) edge (1);
\draw[-stealth,semithick] (2) edge (3);
\draw[-stealth,semithick] (4) edge (2);
\draw[-stealth,semithick] (4) edge (3);
\node[circle, inner sep=1pt, label={$B_4$}] (a) at (1.5,-4.5) {};
\end{tikzpicture}\quad\quad
\begin{tikzpicture}[thick,scale=0.6]
\node[fill,circle,inner sep=0pt,minimum size=3.5pt] (1) at (0,0){}; 
\node[fill,circle,inner sep=0pt,minimum size=3.5pt] (2) at (3,0){}; 
\node[fill,circle,inner sep=0pt,minimum size=3.5pt] (3) at (0,-3){}; 
\node[fill,circle,inner sep=0pt,minimum size=3.5pt] (4) at (3,-3){}; 
\draw[-stealth,semithick] (1) edge (2);
\draw[-stealth,semithick] (3) edge (1);
\draw[-stealth,semithick] (1) edge (4);
\draw[-stealth,semithick] (2) edge (3);
\draw[-stealth,semithick] (4) edge (2);
\draw[-stealth,semithick] (3) edge (4);
\node[circle, inner sep=1pt, label={$C_4$}] (a) at (1.5,-4.5) {};
\end{tikzpicture}\quad\quad
 \begin{tikzpicture}[thick,scale=0.6]
\node[fill,circle,inner sep=0pt,minimum size=3.5pt] (1) at (0,0){}; 
\node[fill,circle,inner sep=0pt,minimum size=3.5pt] (2) at (3,0){}; 
\node[fill,circle,inner sep=0pt,minimum size=3.5pt] (3) at (0,-3){}; 
\node[fill,circle,inner sep=0pt,minimum size=3.5pt] (4) at (3,-3){}; 
\draw[-stealth,semithick] (1) edge (2);
\draw[-stealth,semithick] (3) edge (1);
\draw[-stealth,semithick] (1) edge (4);
\draw[-stealth,semithick] (2) edge (3);
\draw[-stealth,semithick] (2) edge (4);
\draw[-stealth,semithick] (3) edge (4);
\node[circle, inner sep=1pt, label={$D_4$}] (a) at (1.5,-4.5) {};
\end{tikzpicture}\quad\quad\\
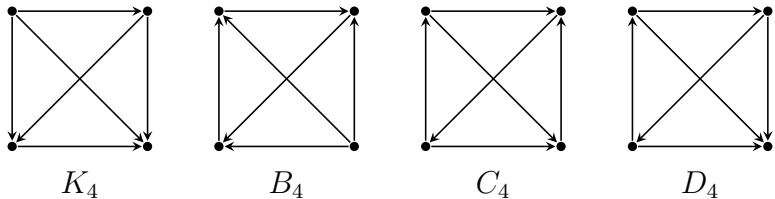
Figure 1: All tournaments with 4 vertices up to isomorphism.
\end{center}
The following result is easy to prove but worth mentioning:
\begin{theorem}\label{THM1.2}
    WMISP for $H$-free tournaments can be solved in polynomial time if $|V(H)| = 4$.
\end{theorem}
\proof{} There are finitely many $K_4$-free tournaments, so WMISP on $K_4$-free tournaments can be solved in constant time. Note that MISP on $H$-free tournaments is equivalent to MISP on $\bar{H}$-free tournaments, since a transitive subtournament is self-complementary. Note $\bar{B_4} \cong D_4$, so it suffices to show WMISP on $B_4$-free tournaments is in P and WMSIP on $C_4$-free tournaments is in P.

Let $G$ be a $B_4$-free tournament that is not transitive with weight assignment $w$, and let $S^* = \{v_1,...,v_k\}$ be an optimal solution of WMISP on $G$ such that $v_i \rightarrow v_j$ if and only if $i< j$. Then, $\{v_2,...,v_k\} \in N^+(v_1)$. Note that, since $G$ is $B_4$-free, $N^+(v)$ is an independent set for all $v \in V(G)$. Thus, we may assume $N^+(v_1) = \{v_2,...,v_k\}$. Hence, $v_1 \in \arg\max_{v \in V(G)} w(N^+(v) \cup \{v\})$. There are only $n$ options for $v_1$, so we can check the weight of all sets in $\{N^+(v) \cup \{v\}:v \in V(G)\}$ to find the optimal solution. It is clear that this can be accomplished in polynomial time. 

Now, let $G$ be a $C_4$-free tournament that is not transitive with weight assignment $w$. Let $T = \{a,b,c\}$ be a cyclic triangle of $G$ with $a \rightarrow b \rightarrow c$. Since $G$ is $C_4$-free, every vertex in $G \setminus T$ is either a common out-neighbor of $T$ or a common in-neighbor of $T$. Let $T^+$ be the common out-neighbors of $T$ and $T^-$ be the common in-neighbor of $T$. We claim that $T^- \Rightarrow T^+$. Suppose otherwise, and let $v_1 \in T^+, v_2 \in T^-$ with $v_1 \rightarrow v_2$, then $\{v_1,v_2,a,b\}$ induces a $C_4$ in $G$. This means we can decompose every $C_4$-free tournaments into a sequence of vertex subsets $S_1,...,S_k$ such that each $S_i$ is either a cyclic triangle or an independent set, and $S_i \Rightarrow S_j$ if and only if $i < j$. In other words, each strong component of $G$ has size 1 or 3. For each $i$ such that $S_i$ induces a cyclic triangle, let $S_i = \{a_i,b_i,c_i\}$ with $w(a_i) \geq w(b_i) \geq w(c_i)$. Let $S = \{S_i: i \in [k], S_i \text{ independent set}\}$ and let $S' = \bigcup _{S_i \text{ cyclic triangle}}\{a_i,b_i\}$. It is clear that $S \cup S'$ is an optimal solution of WMISP on $G$ and $S \cup S'$ can be found in polynomial time. $\qed$

\subsection{Our Main Results}
We are interested in the cases where $H$ has 5 vertices. It is generally more difficult to characterize $H$-free tournaments when $H$ has 5 vertices. However, for some $H$ with 5 vertices, the class of $H$-free prime tournaments has been studied. (We will define prime tournaments in section 2.) For example, Latka gives a structure theorem for $W_5$-free prime tournaments in \cite{Latka} and Liu gives a structure theorem for $U_5$-free prime tournaments. We show that we can reduce WMISP on $H$-free tournaments to WMISP on prime $H$-free tournaments. We further show some complexity results for WMISP on prime $H$-free tournaments with those structure theorems.\\

\noindent Here are the five vertex tournaments $T_5,U_5,W_5$:
\begin{center}
\begin{tikzpicture}[thick,scale=0.6]
\node[fill,circle,inner sep=0pt,minimum size=5pt] (1) at (0,3){}; 
\node[circle,inner sep=0pt,minimum size=5pt] (l1) at (0,3.5){$v_1$}; 
\node[fill,circle,inner sep=0pt,minimum size=5pt] (2) at (2.853,0.927){}; 
\node[circle,inner sep=0pt,minimum size=5pt] (l2) at (3.453,0.927){$v_2$};; 
\node[fill,circle,inner sep=0pt,minimum size=5pt] (3) at (1.763,-2.427){}; 
\node[circle,inner sep=0pt,minimum size=5pt] (l3) at (2.363,-2.427){$v_3$};; 
\node[fill,circle,inner sep=0pt,minimum size=5pt] (4) at (-1.763,-2.427){}; 
\node[circle,inner sep=0pt,minimum size=5pt] (l4) at (-2.363,-2.427){$v_4$};; 
\node[fill,circle,inner sep=0pt,minimum size=5pt] (5) at (-2.853,0.927){}; 
\node[circle,inner sep=0pt,minimum size=5pt] (l5) at (-3.453,0.927){$v_5$};; 
\draw[-stealth,semithick] (1) edge (2);
\draw[-stealth,semithick] (1) edge (3);
\draw[-stealth,semithick] (4) edge (1);
\draw[-stealth,semithick] (5) edge (1);
\draw[-stealth,semithick] (2) edge (3);
\draw[-stealth,semithick] (2) edge (4);
\draw[-stealth,semithick] (5) edge (2);
\draw[-stealth,semithick] (3) edge (4);
\draw[-stealth,semithick] (3) edge (5);
\draw[-stealth,semithick] (4) edge (5);
\node[circle, inner sep=1pt, label={$T_5$}] (a) at (0,-4) {};
\end{tikzpicture}\quad\quad
\begin{tikzpicture}[thick,scale=0.6]
\node[fill,circle,inner sep=0pt,minimum size=5pt] (1) at (0,3){}; 
\node[circle,inner sep=0pt,minimum size=5pt] (l1) at (0,3.5){$v_1$}; 
\node[fill,circle,inner sep=0pt,minimum size=5pt] (2) at (2.853,0.927){}; 
\node[circle,inner sep=0pt,minimum size=5pt] (l2) at (3.453,0.927){$v_2$};; 
\node[fill,circle,inner sep=0pt,minimum size=5pt] (3) at (1.763,-2.427){}; 
\node[circle,inner sep=0pt,minimum size=5pt] (l3) at (2.363,-2.427){$v_3$};; 
\node[fill,circle,inner sep=0pt,minimum size=5pt] (4) at (-1.763,-2.427){}; 
\node[circle,inner sep=0pt,minimum size=5pt] (l4) at (-2.363,-2.427){$v_4$};; 
\node[fill,circle,inner sep=0pt,minimum size=5pt] (5) at (-2.853,0.927){}; 
\node[circle,inner sep=0pt,minimum size=5pt] (l5) at (-3.453,0.927){$v_5$};; 
\draw[-stealth,semithick] (2) edge (1);
\draw[-stealth,semithick] (1) edge (3);
\draw[-stealth,semithick] (4) edge (1);
\draw[-stealth,semithick] (5) edge (1);
\draw[-stealth,semithick] (2) edge (3);
\draw[-stealth,semithick] (2) edge (4);
\draw[-stealth,semithick] (5) edge (2);
\draw[-stealth,semithick] (3) edge (4);
\draw[-stealth,semithick] (3) edge (5);
\draw[-stealth,semithick] (4) edge (5);
\node[circle, inner sep=1pt, label={$U_5$}] (a) at (0,-4) {};
\end{tikzpicture}  \quad\quad
\begin{tikzpicture}[thick,scale=0.6]
\node[fill,circle,inner sep=0pt,minimum size=5pt] (1) at (0,3){}; 
\node[circle,inner sep=0pt,minimum size=5pt] (l1) at (0,3.5){$v$}; 
\node[fill,circle,inner sep=0pt,minimum size=5pt] (2) at (2.853,0.927){}; 
\node[circle,inner sep=0pt,minimum size=5pt] (l2) at (3.453,0.927){$w_4$};; 
\node[fill,circle,inner sep=0pt,minimum size=5pt] (3) at (1.763,-2.427){}; 
\node[circle,inner sep=0pt,minimum size=5pt] (l3) at (2.363,-2.427){$w_3$};; 
\node[fill,circle,inner sep=0pt,minimum size=5pt] (4) at (-1.763,-2.427){}; 
\node[circle,inner sep=0pt,minimum size=5pt] (l4) at (-2.363,-2.427){$w_2$};; 
\node[fill,circle,inner sep=0pt,minimum size=5pt] (5) at (-2.853,0.927){}; 
\node[circle,inner sep=0pt,minimum size=5pt] (l5) at (-3.453,0.927){$w_1$};; 
\draw[-stealth,semithick] (2) edge (1);
\draw[-stealth,semithick] (1) edge (3);
\draw[-stealth,semithick] (4) edge (1);
\draw[-stealth,semithick] (1) edge (5);
\draw[-stealth,semithick] (3) edge (2);
\draw[-stealth,semithick] (4) edge (2);
\draw[-stealth,semithick] (5) edge (2);
\draw[-stealth,semithick] (4) edge (3);
\draw[-stealth,semithick] (5) edge (3);
\draw[-stealth,semithick] (5) edge (4);
\node[circle, inner sep=1pt, label={$W_5$}] (a) at (0,-4) {};
\end{tikzpicture}\\
Figure 2: $T_5,U_5$ and $W_5$.
\end{center}

\noindent In Section \ref{2.2}, we will show:
\begin{theorem}\label{w5f}
WMISP for $W_5$-free tournaments can be solved in polynomial time.
\end{theorem}
\noindent In Section \ref{2.3}, we will show:
\begin{theorem}\label{u5f}
WMISP for $U_5$-free tournaments can be solved in polynomial time.
\end{theorem}

We also give a necessary condition for tournaments $H$ such that MISP restricted to $H$-free tournaments can be solved in polynomial time.

\begin{sdef}
    A tournament $G$ is 1-out(in)-degenerate if every subtournament $G'$ of $G$ has a vertex of out(in)-degree at most $1$.
\end{sdef}
\begin{theorem}\label{1dg}
If MISP restricted to $H$-free tournaments can be solved in polynomial time, then $H$ is both 1-in-degenerate and 1-out-degenerate.
\end{theorem}

Note that $T_5$ is not 1-in-degenerate, so we have the following corollary:
\begin{scorollary}\label{t5f}
MISP for $T_5$-free tournaments is NP-complete.
\end{scorollary}
\newpage
\section{Polynomial-time algorithms}
\subsection{Polynomial-time reduction to prime tournaments}
A \textit{homogeneous set} $X \subseteq V(G)$ is a set of vertices such that for each vertex $v \in V(G) \backslash X$, either $v \Rightarrow X$ or $X \Rightarrow v$. A homogeneous set $X$ is \textit{trivial} if $|X| \leq 1$ or $|X| = |V(G)|$, otherwise it is \textit{nontrivial}.
A tournament is \textit{prime} if it has no nontrivial homogeneous set.

We first show the following theorem via a proof analogous the graph version of the result:
\begin{stheorem}\label{reduce2prime}
If WMISP is in $P$ for prime $H$-free tournaments, then WMISP is in $P$ for all $H$-free tournaments.
\end{stheorem}
The following result is useful to prove \textbf{Theorem 2.1.1}.
\begin{stheorem}[Haglin and Wolf,\cite{MJ96}]\label{findhomoset}
All minimal non-trivial homogeneous sets of a tournament with $n$ vertices can be found in $O(n^4)$.
\end{stheorem}
\proof{ of Theorem 2.1.1}
Suppose we want to solve WMISP for a weighted $H$-free tournament $G$ which is not prime. Let the weight function of $G$ be $w$. Let $X \subset V(G)$ be a smallest nontrivial homogeneous set. We claim that $X$ induces a prime tournament: If there exists $U \subseteq X$ that is a nontrivial homogeneous set in $G[X]$, then $U$ is also a nontrivial homogeneous set in $G$, so $U$ is a smaller nontrivial homogeneous set in $G$, contradicting that $X$ is smallest. Therefore $X$ induces a prime tournament. By Theorem \ref{findhomoset},  such $X$ can be found in $O(n^4)$. 

We claim that we can reduce the problem to WMISP of a tournament smaller than $G$. Let $S$ be an independent subset of $V(G)$ of maximum weight. If $S \cap X$ is non-empty, then we may partition $S$ into $3$ sets $S_1,S_2, S \cap X$ such that $S_1 \Rightarrow S \cap X \Rightarrow S_2$. We claim that $S \cap X$ is a transitive subtournament of $X$ of maximum weight. Suppose for contradiction that there exists a non-empty independent set $U \subseteq X$ such that $w(U) > w(S \cap X)$, then we have $w(S_1 \cup U \cup S_2) > w(S)$ and $S_1 \cup X \cup S_2$ is independent, which contradicts the choice of $S$. Let $w^*(X)$ denote the maximum weight of independent set in $X$. Thus either $S$  intersects $X$ with a set of weight $w^*(X)$ or $S$ is disjoint from $X$. Now we replace the set $X$ with a new single vertex $x$ of weight $w^*(X)$ and such that for all $v \in V(G) \backslash X$, we have $v \rightarrow x$ if and only if $v \Rightarrow X$. We denote the new tournament as $G/X$, with weight function $w'$. If $S$ is not disjoint from $X$, then $S \setminus X \cup \{x\}$ is an independent set of weight $w(S)$ in $G/X$; if $S$ is disjoint from $S$ then $S$ still exists in $G/X$. On the other hand, let $R$ be a maximum weight independent set in $G/X$, if $x \notin R$ then $R$ exists in $G$ with same weight; if $x \in R$ then $(R \setminus \{x\}) \cup (S \cap X)$ is an independent set of same weight in $G$. Thus WMISP of $G/X$ returns the same value as WMISP of $G$. Since $X$ is nontrival, $|X| > 1$, we have $|V(G/X)| < |V(G)|$.

Now, suppose we are given an algorithm A to solve WMISP on prime $H$-free tournaments in polynomial time, then we can do the following on a general tournament: We first find a minimal non-trivial homogeneous set $X$ and compute the optimal solution $S$ of WMISP on $G[X]$ using the algorithm A. Then, we replace $X$ by a single vertex with weight $w(S)$, to get the graph $G/X$. We repeat this process until we reach a prime tournament (i.e., there is no non-trivial homogeneous set). We then solve WMISP on the prime tournament with algorithm A, which gives the optimum solution of WMISP on our original tournament.

We claim that if the running time for solving WMISP on a prime tournament with $n$ vertices is $f(n)$, then the running time for WMISP is at most $O(n^5)+O(n)f(n)$. Notice that this algorithm is recursive, but since $|G/X| \leq  |G| - 1$, the function recurs for at most $n$ times. In each iteration, we work on a tournament with at most $n$ vertices, and in each iteration we take two actions: We find a minimal nontrivial homogeneous set which takes $O(n^4)$ time, and we compute the optimum solution on the tournament $G[X]$ which takes $f(n)$ time. Thus the total running time is bounded by $O(n^5)+O(n)f(n)$. Thus, if $f$ is a polynomial, then the algorithm runs in polynomial time.  $\qed$

\subsection{Excluding $W_5$\label{2.2}}
Following the notation of \cite{Liu_U5}, we define the following tournaments:
\begin{sdef}
    Let $n \geq 1$ be odd. The tournaments $T_n, U_n, W_n$ are defined as:
    \begin{itemize}
    \item 
    $T_n$ is the tournament with vertex set $\{v_1,...,v_n\}$ such that $v_i \rightarrow v_j$ if $j - i \equiv 1,2,...,\frac{n-1}{2} (\text{mod } n)$.
    \item 
    $U_n$ is the tournament obtained from $T_n$ by reversing all edges with both ends in $\{v_1,v_2,...,v_{(n-1)/2}\}$.
    \item 
    $W_n$ is the tournament with vertex set $\{v,w_1,...,w_{n-1}\}$ such that $w_i \rightarrow w_j$ if $i < j$, and $\{w_i: i\ even\} \Rightarrow v \Rightarrow \{w_i: i\ odd\}$.
    \end{itemize}
\end{sdef}
\noindent Now we define a few other tournaments. Let $I_n$ be the transitive tournament on $n$ vertices. Let $Q_7$ be the tournament with vertex set $\mathbb{Z}_7$ and edge set $\{(i,j):j-i \in \{1,2,4\}\}$. Since $Q_7$ is vertex-transitive, we can define $Q_7-v$ as the tournament obtained from deleting a vertex of $Q_7$. $Q_7$ is also known as the 7-vertex Paley tournament.
\begin{center}
\begin{tikzpicture}[thick,scale=0.6]
\node[fill,circle,inner sep=0pt,minimum size=5pt] (1) at (0,3){}; 
\node[circle,inner sep=0pt,minimum size=5pt] (l1) at (0,3.5){$v_0$}; 
\node[fill,circle,inner sep=0pt,minimum size=5pt] (2) at (3.063,1.524){}; 
\node[circle,inner sep=0pt,minimum size=5pt] (l2) at (3.593,1.524){$v_1$}; 
\node[fill,circle,inner sep=0pt,minimum size=5pt] (3) at (3.820,-1.789){}; 
\node[circle,inner sep=0pt,minimum size=5pt] (l3) at (4.370,-1.789){$v_2$}; 
\node[fill,circle,inner sep=0pt,minimum size=5pt] (4) at (1.700,-4.449){}; 
\node[circle,inner sep=0pt,minimum size=5pt] (l4) at (2.300,-4.549){$v_3$}; 
\node[fill,circle,inner sep=0pt,minimum size=5pt] (5) at (-1.700,-4.449){}; 
\node[circle,inner sep=0pt,minimum size=5pt] (l5) at (-2.300,-4.549){$v_4$};
\node[fill,circle,inner sep=0pt,minimum size=5pt] (6) at (-3.820,-1.789){}; 
\node[circle,inner sep=0pt,minimum size=5pt] (l6) at (-4.370,-1.789){$v_5$}; 
\node[fill,circle,inner sep=0pt,minimum size=5pt] (7) at (-3.063,1.524){}; 
\node[circle,inner sep=0pt,minimum size=5pt] (l7) at (-3.593,1.524){$v_6$}; 
\draw[-stealth,semithick] (1) edge (2);
\draw[-stealth,semithick] (2) edge (3);
\draw[-stealth,semithick] (3) edge (4);
\draw[-stealth,semithick] (4) edge (5);
\draw[-stealth,semithick] (5) edge (6);
\draw[-stealth,semithick] (6) edge (7);
\draw[-stealth,semithick] (7) edge (1);
\draw[-stealth,semithick] (1) edge (3);
\draw[-stealth,semithick] (2) edge (4);
\draw[-stealth,semithick] (3) edge (5);
\draw[-stealth,semithick] (4) edge (6);
\draw[-stealth,semithick] (5) edge (7);
\draw[-stealth,semithick] (6) edge (1);
\draw[-stealth,semithick] (7) edge (2);
\draw[-stealth,semithick] (1) edge (5);
\draw[-stealth,semithick] (2) edge (6);
\draw[-stealth,semithick] (3) edge (7);
\draw[-stealth,semithick] (4) edge (1);
\draw[-stealth,semithick] (5) edge (2);
\draw[-stealth,semithick] (6) edge (3);
\draw[-stealth,semithick] (7) edge (4);
\node[circle, inner sep=1pt, label={$Q_7$}] (a) at (0,-6) {};
\end{tikzpicture}\qquad\qquad\qquad\quad
\begin{tikzpicture}[thick,scale=0.6]
\node[fill,circle,inner sep=0pt,minimum size=5pt] (1) at (0,3){}; 
\node[circle,inner sep=0pt,minimum size=5pt] (l1) at (0,3.5){$v_0$}; 
\node[fill,circle,inner sep=0pt,minimum size=5pt] (2) at (2.944,1.300){}; 
\node[circle,inner sep=0pt,minimum size=5pt] (l2) at (3.474,1.300){$v_1$}; 
\node[fill,circle,inner sep=0pt,minimum size=5pt] (3) at (2.944,-2.100){}; 
\node[circle,inner sep=0pt,minimum size=5pt] (l2) at (3.474,-2.100){$v_2$}; 
\node[fill,circle,inner sep=0pt,minimum size=5pt] (4) at (0.000,-3.800){}; 
\node[circle,inner sep=0pt,minimum size=5pt] (l1) at (0,-4.3){$v_3$}; 
\node[fill,circle,inner sep=0pt,minimum size=5pt] (5) at (-2.944,-2.100){};
\node[circle,inner sep=0pt,minimum size=5pt] (l2) at (-3.474,-2.100){$v_4$}; 
\node[fill,circle,inner sep=0pt,minimum size=5pt] (6) at (-2.944,1.300){}; 
\node[circle,inner sep=0pt,minimum size=5pt] (l6) at (-3.474,1.300){$v_5$}; 
\draw[-stealth,semithick] (1) edge (2);
\draw[-stealth,semithick] (2) edge (3);
\draw[-stealth,semithick] (3) edge (4);
\draw[-stealth,semithick] (4) edge (5);
\draw[-stealth,semithick] (5) edge (6);

\draw[-stealth,semithick] (1) edge (3);
\draw[-stealth,semithick] (2) edge (4);
\draw[-stealth,semithick] (3) edge (5);
\draw[-stealth,semithick] (4) edge (6);

\draw[-stealth,semithick] (6) edge (1);

\draw[-stealth,semithick] (1) edge (5);
\draw[-stealth,semithick] (2) edge (6);

\draw[-stealth,semithick] (4) edge (1);
\draw[-stealth,semithick] (5) edge (2);
\draw[-stealth,semithick] (6) edge (3);

\node[circle, inner sep=1pt, label={$Q_7-v$}] (a) at (0,-6) {};
\end{tikzpicture} \\
Figure 3: $Q_7$ and $Q_7-v$
\end{center}
\begin{stheorem}[Latka, \cite{Latka}]\label{Structure-w5f}
A prime tournament is $W_5$-free if and only if it is isomorphic to one of $I_1,I_2,Q_7-v, Q_7, T_n$ or $U_n$ for some odd $n \geq 1$.
\end{stheorem}
We now introduce two lemmas in order to prove \textbf{Theorem 1.2}.
\begin{slemma}\label{Tn->P}
WMISP for the class $\{T_n: n\geq 1, n\text{ odd}\}$ can be solved in polynomial time.
\end{slemma}
\proof{} In \cite{LIU-VariousThmsonTounrmanets}, Liu showed that $T_n$ is $D_4$-free. Then, by \ref{THM1.2}, WMISP for  $\{T_n: n\geq 1, n\text{ odd}\}$ can be solved in polynomial time. $\qed$
\begin{slemma}\label{Un->P}
WMISP for the class $\{U_n: n\geq 1, n\text{ odd}\}$ can be solved in polynomial time.
\end{slemma}
\proof{} Let $n \geq 1$ be an odd number and suppose we want to solve WMISP of $G \cong U_n$, with a weight assignment $w$ on $V(G)$. We may assume $n \geq 5$. We first aim to label the vertices as $v_1,...,v_n$ following the definition of $U_n$. For short, we call this a proper labeling of $V(G)$. Here are two useful properties from the definition of $U_n$ for $n \geq 5$:
\begin{itemize}
\item $v_{\frac{n-1}{2}}$ is the unique vertex with in-degree $1$, and its in-neighbor is $v_n$;
\item $U_n \setminus \{v_{\frac{n-1}{2}},v_n\} \cong U_{n-2}$ for $n \geq 5$.
\end{itemize}
It is then clear that we can find a correct way to label the graph $G$ in polynomial time.

Now with the proper labeling, we let $A = \{v_{\frac{n+1}{2}},...,v_n\}$ and $B = \{v_1,...,v_{\frac{n-1}{2}}\}$. Let $S$ be a maximum-weight independent set of $G \cong U_n$. If $S \subseteq A$ or $S \subseteq B$, then since $A,B$ are both independent, we have $S = A$ or $S = B$. Otherwise,  we may assume $S \cap A \neq \varnothing$ and $S \cap B \neq \varnothing$. \\

We first claim that $U_n$ has no cyclic triangle with exactly two vertices in $B$. Suppose $v_i,v_j,v_k$ induces a cyclic triangle and $v_i \in A, v_j,v_k \in B$ with $j > k$. Since $v_i \in A$, $i > j > k$. Then $v_j \rightarrow v_k$ since $v_j,v_k \in B$ with $j > k$, and therefore $v_k \rightarrow v_i, v_i \rightarrow v_j$, which means $i - k < \frac{n}{2}$ and $i - j > \frac{n}{2}$. This contradicts $j > k$.\\

If $|S \cap A| = 1$, then since $U_n$ has no cyclic triangle with exactly two vertices in $B$, it follows that $S \cap B = B$. Now assume that $|S \cap A| \geq 2$. Let $v_i$ be the vertex of smallest index in $S \cap A$ and $v_j$ be the vertex of largest index in $S \cap A$. We claim that $S' = \{v_i,v_{i+1},...,v_j\} \subseteq S$. Suppose $v_k \notin S$ and $i<k<j$, then by maximality of $S$, it follows that $v_k$ induces a cyclic triangle with two vertices in $S$. Let $v_k,v_p,v_q$ be a cyclic triangle such that $v_p,v_q \in S$. Since $A$ is acyclic, we have $\{v_p,v_q\} \cap B \neq \varnothing$. Without loss of generality, we may assume $p < q$. Then either $\{v_p,v_q\} \subseteq B$ or $v_p \in B,v_q \in A$. We consider three cases here. First assume that $v_p,v_q \in B$. This is impossible since there is no cyclic triangle with two vertices in $B$. Second, assume that $v_p \in B, v_q \in A$ and $q < k$. Then $v_q \rightarrow v_k, v_k \rightarrow v_p, v_p \rightarrow v_q$. This implies $k-p > \frac{n}{2}$, and therefore $j - p > \frac{n}{2}$. As a result, $\{v_j,v_p,v_q\}$ is a cyclic triangle, which contradicts that $S$ is independent. Third, assume that $v_p \in B, v_q \in A$ and $q > k$, then a similar argument will give that $\{v_i,v_p,v_q\}$ is a cyclic triangle. Thus such $v_k$ does not exist, which means $S' = \{v_i,v_{i+1},...,v_j\} \subseteq S$. This suggests that we only need to find $i,j$ to fix $S \cap A$.\\

Now suppose we fix $i,j$ (and therefore fix $S \cap A$). We claim that $S \cap B= B\setminus \{v_k:v_j \rightarrow v_k,v_k \rightarrow v_i\}$. If $v_j \rightarrow v_k$ and $v_k \rightarrow v_i$, then $\{v_j,v_k,v_i\}$ induces a cyclic triangle, and thus $v_k \notin S$. On the other hand, suppose $v_l \in B\setminus \{v_k:v_j \rightarrow v_k,v_k \rightarrow v_i\}$ and $v_k \notin S$. Then since no cyclic triangle has two vertices in $B$, it follows that $v_l$ induces a cyclic triangle with two vertices in $S \cap A$. Let $\{v_p,v_q,v_l\}$ be such a triangle with $i \leq p < q \leq j$, then $\{v_i,v_j,v_l\}$ is also a cyclic triangle, so $v_l \in  \{v_k:v_j \rightarrow v_k,v_k \rightarrow v_i\}$, which is a contradiction. Hence $S \cap B= B\setminus \{v_k:v_j \rightarrow v_k,v_k \rightarrow v_i\}$. Note that by definition of $U_n$ we have $B\setminus \{v_k:v_j \rightarrow v_k,v_k \rightarrow v_i\}= B\setminus \{v_k:i-\frac{n}{2} < k < j-\frac{n}{2}\}$.\\

Thus in total we have 3 cases for $S$:
\begin{enumerate}
    \item  $S = A$ or $S = B$
    \item  $|S \cap A| = 1$ and $S \cap B = B$
    \item  $S \cap A = \{v_i,v_{i+1},...,v_j\}$ for some $\frac{n+1}{2} \leq i < j \leq n $ and $S \cap B = B\setminus\{v_k:i-\frac{n}{2} < k < j-\frac{n}{2}\}$
\end{enumerate}
There are $O(1)$ sets to check for Case $1$, $O(n)$ sets to check for Case 2 and $O(n^2)$ sets to check for Case 3. It is obvious that the weight of all these sets can be computed in polynomial time with the proper labelling of $U_n$ given, and we have showed that we can compute the proper labelling in polynomial time. Thus we can enumerate all possible candidates for $S$ in polynomial time. Thus, compute the maximum weight of an independent set in $U_n$ is in polynomial time. $\qed$\\

Now we can prove \textbf{Theorem \ref{w5f}}:
\begin{customthm}{\ref{w5f}}\label{THM1.2}
WMISP for $W_5$-free tournaments can be solved in polynomial time.
\end{customthm}

\proof{}
By Theorem \ref{reduce2prime}, it suffices to show that we can solve WMISP in polynomial time for prime $W_5$-free tournaments. Let $G$ be a $W_5$-free transitive tournament with $n \geq 8$ vertices. By Theorem \ref{Structure-w5f}, $G$ is isomorphic to $T_n$ or $U_n$ for some odd $n \geq 1$. We can check whether all vertices have same out-degree in polynomial time. If all vertices have same out-degree, then $G \cong T_n$, and by Lemma \ref{Tn->P}, we can solve WMISP on $G$ in polynomial time. Otherwise, $G \cong U_n$, and by Lemma \ref{Un->P}, we can solve WMISP on $G$ in polynomial time. Thus MISP for $W_5$-free tournaments can be solved in polynomial time. $\qed$

\subsection{Excluding $U_5$\label{2.3}}
Here is the structure theorem for all prime $U_5$-free tournaments:
\begin{stheorem}[Liu, \cite{Liu_U5}]\label{Structure_u5f}
Let G be a prime tournament. Then $G$ is $U_5$-free if and only if $G$ is $T_n$ for some odd $n \geq 1$ or $V(G)$ can be partitioned into sets $ X, Y, Z $ such that $ X \cup Y , Y \cup Z, \text{and } Z \cup X$ are all independent.
\end{stheorem}
\begin{slemma}\label{FindXYZ}
There is a polynomial algorithm that takes a prime $U_5$-free tournament $G$ as input, and returns a partition $(X,Y,Z)$ of $V(G)$ such that $ X \cup Y , Y \cup Z, \text{and } Z \cup X$ are all independent.    
\end{slemma}
\proof{}
 The following algorithm is due to Seokbeom Kim (private communication):
\begin{algorithm}[H]
\caption*{\textbf{$(X,Y,Z)$-partition}}
\begin{algorithmic}
    \Require $G \ $is a $\ U_5$-free, prime tournament and $G[T]$ is a cyclic triangle
    \Function{FindPartition}{G,T}
     \State{$T = \{x_0,y_0,z_0\}$, $X = \{x_0\}, Y = \{y_0\}, Z = \{z_0\}$, $S = X \cup Y \cup Z$}
     \While{$S \neq V(G)$}
        \For{$i \in V(G) \setminus S$}
            \If{$i$ induces a cyclic triangle with $S$}
                \If{ $X \cup Y \cup \{i\}$ contains a cyclic triangle} \State{$Z \gets Z \cup \{i\}$}
                \ElsIf{$X \cup Z \cup \{i\}$ contains a cyclic triangle} \State{$Y \gets Y \cup \{i\}$}
                \Else \State{$X \gets X \cup \{i\}$}
                \EndIf
            \State{$S \gets S \cup \{i\}$}
            \EndIf
        \EndFor
    \EndWhile\\
    \textbf{return} $(X,Y,Z)$
    \EndFunction
\end{algorithmic}
\end{algorithm}
We claim this algorithm will terminate and return a correct partition for the structure theorem. We notice the label of the three vertices in $T$ must be correct up to symmetry. Now, inductively at some step of the algorithm, suppose all vertices in $S$ are labeled correctly. If $S \neq V(G)$, then we claim that there is a vertex not in $S$ that induces a cyclic triangle with $S$. Since we start with a cyclic triangle, inductively we assume that $G[S]$ is strongly connected. Since $S\neq V(G)$ and $G$ is prime, it follows that there is a vertex $i \not\in S$ with both an in-neighbor and an out-neighbor in $S$. Since $G[S]$ is strongly connected, and $N^+(v) \cap S, N^-(v) \cap S$ is a partition of $S$ into two non-empty parts, it follows that there exists $u \in N^+(v) \cap S, v \in N^-(v) \cap S$ such that $u \rightarrow v$. Now this gives a cyclic triangle $\{u,v,i\}$ in $S \cup \{i\}$. Now in $G[S \cup \{i\}]$, the label of all vertices are determined correctly, since exactly one of $u,v,i$ is in $X,Y,Z$ respectively, and we already have $u,v$ labeled correctly. Moreover, $G[S \cup \{i\}]$ is strongly connected, maintaining our inductive assumption. Inductively, this algorithm correctly labels every vertex as one of $X,Y,Z$. It is clear that this algorithm runs in polynomial time. $\qed$\\

The algorithm shows that if we want to determine the complexity of WMISP for $U_5$-free prime tournaments, it is safe to assume that we have the partition $(X,Y,Z)$ as stated in the structure theorem. With this lemma, we are able to show the following result:
\begin{customthm}{\ref{u5f}}\label{THM1.3}
WMISP for $U_5$-free tournaments can be solved in polynomial time.
\end{customthm}

To prove \textbf{Theorem \ref{u5f}}, we first define some terminology that will be used in the proof:
\begin{sdef}
     Let $G$ be a tournament with a tri-partition $(X,Y,Z)$ of $V(G)$. A cyclic triangle is called an $XYZ$ triangle if it consists of vertices $x \in X, y \in Y, z \in Z$ such that $x \rightarrow y, y \rightarrow z, z \rightarrow x$.
\end{sdef}

\begin{slemma}\label{SingleTypeTriangle}
Let G be a tournament. If $v \in V(G)$ and $V(G)$ can be partitioned into sets $ X, Y, Z $ such that $ X \cup Y , Y \cup Z, \text{and } Z \cup X$ are all independent, then either every cyclic triangle containing $v$ is an $XYZ$ triangle or every cyclic triangle containing $v$ is an $XZY$ triangle.
\end{slemma}
\proof{}
Since $X \cup Y, Y \cup Z, Z \cup X$ are all independent sets, all cyclic triangles containing $v$ must contain one vertex from each of $X,Y,Z$. Now, suppose $v \in X$ and there exists $y_1,y_2 \in Y,z_1,z_2 \in Z$ such that $\{v,y_1,z_1\}$ is an $XYZ$ triangle and $\{v,z_2,y_2\}$ is an $XZY$ triangle. Since $v \rightarrow y_1, y_2 \rightarrow v$ and $X \cup Y$ is transitive, we must have $y_2 \rightarrow y_1$. Similarly, $z_1 \rightarrow z_2$. But then notice $z_1 \rightarrow z_2, z_2 \rightarrow y_2, y_2 \rightarrow y_1, y_1 \rightarrow z_1$, so $Y \cup Z$ is not transitive, which is a contradiction.$\qed$
\begin{slemma}\label{XYZ-only}
If $G$ is a prime $U_5$-free tournament, then in polynomial time, we can find a tri-partition $(X,Y,Z)$ of $V(G)$ such that $X \cup Y,Y \cup Z,Z \cup X$ are all independent, and all cyclic triangles in $G$ are all $XYZ$ triangles.
\end{slemma}
\proof{}
In the proof of Lemma \ref{FindXYZ}, we notice that we start with a XYZ triangle $\{x_0,y_0,z_0\}$. We will inductively show that each vertex of $S$ is contained in an $XYZ$ triangle; then the result follows from Lemma \ref{SingleTypeTriangle}. This is true at the start of the algorithm. When we add a vertex $i$ to $S$, it follows that there are vertices $p,q \in S$ which form a cyclic triangle with $i$. By Lemma \ref{SingleTypeTriangle}, since $p \in S$, it follows that this triangle is an $XYZ$ triangle. $\qed$\\

Therefore, from now on, we may assume $G$ is a prime $U_5$-free tournament with partition $(X,Y,Z)$ such that all cyclic triangles in $G$ are $XYZ$. Now, let $X = \{x_0,x_1,...,x_p\}, Y = \{y_1,...,y_q\}, Z = \{z_1,...,z_r\}$ such that the indices provide a natural topological ordering of each set. We will describe a polynomial-time algorithm that returns the maximum-weight independent set in $G$ using dynamic programming. We first prove a lemma which is useful for our algorithm:
\begin{slemma}\label{Extend-One-Vx}
    Let $G,X,Y,Z$ be defined as above. Let $S_{i,j,k} = \{x_i,x_{i+1},...,x_p\} \cup \{y_j,y_{j+1},...,y_q\} \cup \{z_k,z_{k+1},...,z_r\}$. Given optimum WMISP solutions for all $G[S_{i,j,k}]$ where $i \in [p], j \in [q], k \in [r]$, we can compute optimal WMISP solution for $G$ in polynomial time.
\end{slemma}
\proof{}
Let $S^*$ be an optimum solution of WMISP on $G$. We consider the following 3 cases for computing the optimum solution:
\begin{itemize}
    \item \textbf{Case 1:} $x_0 \notin S^*$
    \item \textbf{Case 2:} $x_0 \in S^*, N^+(x_0) \cap Y \cap S^* \neq \varnothing$
    \item \textbf{Case 3:} $x_0 \in S^*, N^+(x_0) \cap Y \cap S^* = \varnothing$
\end{itemize}
\noindent \textbf{Case 1:} $x_0 \notin S^*$:\\
Then $S^*$ is also an optimum solution of $G[S_{1,1,1}]$. We can simply check the optimum solution of $G[S_{1,1,1}]$.

\noindent \textbf{Case 2:} $x_0 \in S^*$ and $ N^+(x_0) \cap Y \cap S^* \neq \varnothing$:\\
Let $Y'$ be the set of in-neighbours of $x_0$ in $Y$ and let $m$ be the minimum index such that $y_m \in N^+(x_0) \cap S^*$. Let $Z'$ be the set of all in-neighbours of $y_m$. Let $A = X \setminus \{x_0\}$, $B = \{y_{m},...,y_q\}$, and let $C$ be the set of outneighbours of $y_m$ in $Z$ that are also outneighbours of $x_0$. If $S'$ is an optimum solution of $G[A\cup B \cup C]$, we claim that $R^* = \{x_0\} \cup S' \cup Y' \cup Z'$ is an optimum solution of $G$.\\
\begin{center}
\begin{tikzpicture}[thick,scale=0.6]
\node[fill,circle,inner sep=0pt,minimum size=3.5pt,label = {$x_0$}] (x0) at (-7,3){}; 
\node[fill,circle,inner sep=0pt,minimum size=3.5pt,label = {$x_1$}] (x1) at (2,3){}; 
\node[circle,inner sep=0pt,minimum size=3.5pt,label = {$......$}] (x) at (6,2.5){}; 
\node[fill,circle,inner sep=0pt,minimum size=3.5pt,label = {$x_p$}] (xp) at (10,3){}; 
\draw[decoration={brace,mirror,raise=5pt},decorate] (2,3) -- node[below=6pt] {$A$} (10,3);
\node[fill,circle,inner sep=0pt,minimum size=3.5pt,label = {$y_1$}] (y1) at (0,0){}; 
\node[circle,inner sep=0pt,minimum size=3.5pt,label = {$......$}] (y2) at (2,-0.5){}; 
\node[fill,circle,inner sep=0pt,minimum size=3.5pt,label = {west:$y_s$}] (ys) at (4,0){}; 
\node[fill,circle,inner sep=0pt,minimum size=3.5pt,label = {south:$y_{s+1}$}] (ys1) at (6,0){}; 
\node[circle,inner sep=0pt,minimum size=3.5pt,label = {$......$}] (y) at (8,-0.5){};
\node[fill,circle,inner sep=0pt,minimum size=3.5pt,label = {$y_{m}$}] (ym) at (10,0){};
\node[circle,inner sep=0pt,minimum size=3.5pt,label = {$......$}] (y) at (12,-0.5){};
\node[fill,circle,inner sep=0pt,minimum size=3.5pt,label = {$y_{q}$}] (yq) at (14,0){};
\draw[middlearrow={stealth},semithick,bend right = 10] (y1) to (x0);
\draw[middlearrow={stealth},semithick,bend right = 10] (ys) to (x0);
\draw[middlearrow={stealth},semithick,bend left = 10] (x0) to (ys1);
\draw[middlearrow={stealth},semithick,bend left = 10] (x0) to (ym);
\draw[decoration={brace,mirror,raise=5pt},decorate] (0,0) -- node[below=6pt] {$Y'$} (4,0);
\draw[decoration={brace,mirror,raise=5pt},decorate] (10,0) -- node[below=6pt] {$B$} (14,0);
\node[fill,circle,inner sep=0pt,minimum size=3.5pt,label = {south:$z_1$}] (z1) at (0,-3){}; 
\node[circle,inner sep=0pt,minimum size=3.5pt,label = {$......$}] (x) at (2,-3.5){}; 
\node[fill,circle,inner sep=0pt,minimum size=3.5pt,label = {south:$z_t$}] (zt) at (4,-3){}; 
\node[fill,circle,inner sep=0pt,minimum size=3.5pt,label = {south:$z_{t+1}$}] (zt1) at (6,-3){}; 
\node[circle,inner sep=0pt,minimum size=3.5pt,label = {$......$}] (x) at (8,-3.5){}; 
\node[fill,circle,inner sep=0pt,minimum size=3.5pt,label = {south:$z_u$}] (zu) at (10,-3){}; 
\node[fill,circle,inner sep=0pt,minimum size=3.5pt,label = {south:$z_{u+1}$}] (zu1) at (12,-3){}; 
\node[circle,inner sep=0pt,minimum size=3.5pt,label = {$......$}] (x) at (14,-3.5){}; 
\node[fill,circle,inner sep=0pt,minimum size=3.5pt,label = {south:$z_{r}$}] (zr) at (16,-3){}; 
\draw[middlearrow={stealth},semithick] (z1) to (ym);
\draw[middlearrow={stealth},semithick] (zt) to (ym);
\draw[middlearrow={stealth},semithick] (ym) to (zt1);
\draw[middlearrow={stealth},semithick] (ym) to (zu);
\draw[middlearrow={stealth},semithick] (ym) to (zu1);
\draw[decoration={brace,mirror,raise=5pt},decorate] (0,-3.5) -- node[below=6pt] {$Z'$} (4,-3.5);
\draw[decoration={brace,mirror,raise=5pt},decorate] (12,-3.5) -- node[below=6pt] {$C$} (16,-3.5);
\draw[middlearrow={stealth},semithick,bend left = 35] (z1) to (x0); 
\draw[middlearrow={stealth},semithick,bend left = 30] (zt) to (x0);
\draw[middlearrow={stealth},semithick,bend left = 25] (zt1) to (x0);
\draw[middlearrow={stealth},semithick,bend left = 15] (zu) to (x0);
\draw[middlearrow={stealth},semithick,bend right = 12.5] (x0) to (zu1);
\end{tikzpicture}\\
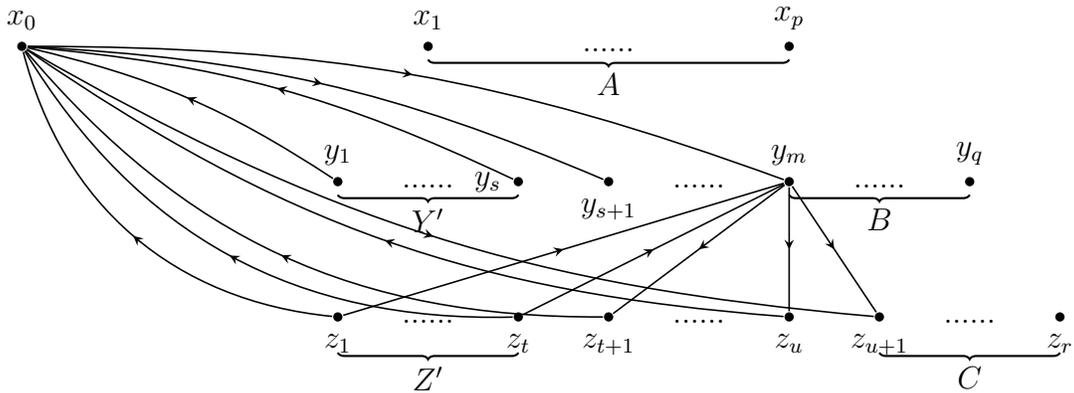
Figure 4: Diagram for Case 2
\end{center}
We first show that $S^* \subseteq  \{x_0\} \cup Y' \cup Z' \cup A \cup B \cup C$. Suppose not; since $A \cup \{x_0\} = X$, we only need to consider $S^* \cap Y$ and $S^* \cap Z$. Since $Y \setminus (Y' \cup B)$ is the set of outneighbours $y_i$ of $x_0$ with index $i < m$, by the definition of $y_m$, we know that $Y \setminus (Y' \cup B)$  is disjoint from $S^*$. For $S^* \cap Z$, we have $Z \setminus (C \cup Z') = Z \cap N^+(y_m) \cap N^-(x_0)$. Since every vertex in this set form a cyclic triangle with $y_m$ and $x_0$, this set is also disjoint from $S^*$, as desired.

In order to prove the claim, it remains to show that $R^*$ is independent in $G$. Since $Y' \Rightarrow \{x_0\} \cup A$ and all triangles are XYZ, there is no cyclic triangle in $R^*$ involving a vertex in $Y'$. Hence, every cyclic triangle in $R^*$ uses a vertex in $B$. Since $Z' \Rightarrow B$, there is no cyclic triangle in $R^*$ using $Z'$. Now, every cyclic triangle in $R^*$ uses a vertex in $C$, and since $x_0 \Rightarrow C$ there is no cyclic triangle in $R^*$ using $x_0$. Hence every cyclic triangle in $R^*$ is in $S'$, and by definition $S'$ is independent. Thus $R^*$ is independent. \\

\noindent\textbf{Case 3}: $x_0 \in S^*, N^+(x_0) \cap Y \cap S^* = \varnothing$:\\
Let $Y'$ be the set of all in-neighbours of $x_0$ in $Y$, we claim that $S^* = X \cup Y' \cup Z$ in this case. Note that since $Y'\Rightarrow x_0 \Rightarrow X \setminus \{x_0\}$ and $X \cup Y$ is independent, it follows that $Y' \Rightarrow X$. Since all triangles are XYZ, it follows that $S^*$ has no cyclic triangles, thus $S^*$ is independent. Notice that $S^*$ contains all vertices that could be in $S^*$ (since we assume $x_0 \in S^*$ and $S^*$ has no outneighbors in $S^* \cap Y$), therefore it is an optimum solution. \\

Now, let $n = |V(G)|$. Note that in Case 1, we only need to check a solution we already have, thus if Case 1 happens, the running time is $O(1)$. For Case 2,  we need to compute the sets $S',Y',Z'$. Consider at most $n$ possible values of $m$; for each, we can ask for $S'$ in $O(1)$ by our assumption, and for the remaining two sets we need $O(n^2)$ times to find them. Hence this requires $O(n)$ running time. For Case 3, we only need to find the set $Y'$ which is also in $O(n)$. Therefore the total running times is polynomial, which completes the proof.$\qed$\\

If we aim to solve WMISP on $G[\{x_i,...,x_p\} \cup \{y_j,...,y_q\} \cup \{z_k,...,z_r\}]$ in polynomial time, we only need to know optimum solutions on $G[\{x_{i'},...,x_p\} \cup \{y_{j'},...,y_q\} \cup \{z_{k'},...,z_r\}]$ for all $i' > i,j' \geq j$ and $k' \geq k$. (By symmetry, this could also be $i' \geq i,j' > j$ and $k' \geq k$ or $i' \geq i,j' \geq j$ and $k' > k$.) Now we are ready to prove the main theorem.\\

\proof{ of \textbf{Theorem} \ref{u5f}}
By Lemma \ref{reduce2prime}, it suffice to show WMISP on prime $U_5$-free tournaments is in P. Given a WMISP on a prime $U_5$-free tournament $G$, we first partition $V(G)$ into three sets $X,Y,Z$ such that $X \cup Y$, $Y\cup Z$ and $Z \cup  X$ are all independent. By Lemma \ref{FindXYZ}, we can find such partition in polynomial time. By Lemma \ref{SingleTypeTriangle}, we may assume all triangles in $G$ are XYZ triangles. Let $X = \{x_1,...,x_p\},Y = \{y_1,...,y_q\}$ and $Z = \{z_1,...,z_r\}$ such that the indices provides a natural topological ordering on each set. For $i \in [1,p], j \in [1,q],k \in [1,r]$, we define $X_i = \{x_i,...,x_p\}, Y_j = \{y_j,...,y_q\},Z_k = \{z_k,...,z_r\}$. We also define $X_{p+1} = Y_{q+1} = Z_{r+1} = \varnothing$.\\

Here is the algorithm:
\begin{algorithm}[H]
\caption*{\textbf{WMISP of prime $U_5$-free tournaments}}
\begin{algorithmic}
\Require{G is prime $U_5$-free with partition $(X,Y,Z)$ as in \ref{Structure_u5f}, with weight function $w:V(G) \rightarrow \mathbb{Q}_+$}
\State{OPT$(a,b,c)$ $\leftarrow 0$ for all $a \in [p],b \in [q], c \in [r]$}
\State{OPT$(a,b,c)$ $\leftarrow w(Y_b)+w(Z_c)$ for all $a = p+1,b \in [q], c \in [r]$}
        \For{$i=p,\ldots,1$}
            \For{$j=q,\ldots,1$}
                \For{$k=r,\ldots,1$}
                    \State{$S_1 = \text{OPT}(i+1,j,k)$}
                    \State{$S_2 = w(X_i \cup Z_k \cup (Y_j \cap N^-(x_i)))$}
                    \For{$y_m \in N^+(x_i) \cap Y_j$}
                        \State{$Y' = N^-(x_i) \cap Y_j$, $Z' = N^-(y_m) \cap Z_k$}
                        \State{$C = Z_k \cap N^+(y_m) \cap N^+(x_i)$}
                        \State{$n \leftarrow$ smallest index of vertices in $C$}
                        \State{$S_{3,m} = w(\{x_i\} \cup Y' \cup Z' \cup \text{OPT}(i+1,m,n))$}
                        \State{$\text{OPT}(i,j,k) = \max\{S_1,S_2,S_{3,m}: m \in [q] \ \ s.t.\ \  y_m \in N^+(x_i) \cap Y_j\}$}
                    \EndFor
                \EndFor
            \EndFor
        \EndFor\\
    \Return OPT$(1,1,1)$
\end{algorithmic}
\end{algorithm}
The correctness of the algorithm above is shown in the proof of Lemma \ref{Extend-One-Vx}. For the running time of the algorithm, we note that we have 4 for-loops, each traverse a subset of vertex, so we have $O(n^4)$ loops to run. For each step inside the loop, note that we need $O(n)$ time to compute $S_{3,m}$ for each $y_m$, and all other steps are in $O(1)$ time. This gives a total running time of $O(n^5)$, which is in polynomial time. $\qed$

\newpage
\section{NP-completeness Results}
In Section 3, we discuss the class of $H$-free tournaments where MISP is still NP-complete.
\subsection{Excluding 1-degenerate Tournaments}
\begin{sdef}
    For a graph $G$ with $n$ vertices, we denote the graph $G\cup nK_1$ (i.e. $G$ with $n$ additional isolated vertices) by $G^+$. 
\end{sdef}
\begin{sdef}
    For a class of graph $\mathcal{C}$, we denote  the class of graphs $\{G^+: G \in \mathcal{C}\}$ by $\mathcal{C}^+$.
\end{sdef}
The following is a corollary from the proof of Theorem \ref{npc}:
\begin{stheorem}\label{VC->MISP}
    If $\mathcal{C}$ is a class of graphs such that MVCP (Minimum Vertex Cover Problem) on $\mathcal{C}$ is NP-complete, then MISP is NP-complete on tournaments with a backedge graph in $\mathcal{C}^+$.
\end{stheorem}
\proof{} This follows directly by using graphs from $\mathcal{C}$ for the reduction in Theorem \ref{npc}. $\qed$\\

We are now ready to prove Theorem \ref{1dg}. The key idea is to use the same reduction used in the proof of Theorem \ref{npc}, and apply the following lemma:
\begin{lemma}[Chv\'atal, Poljak \cite{CHVATAL1973305,Poljak1974}]\label{path3MVCP}
    The minimum vertex cover problem on a graph $G$ is equivalent to the minimum vertex cover problem on $G'$ where $G'$ is obtained by replacing each edge of $G$ by a path of 3 edges.
\end{lemma}

\begin{customthm}{\ref{1dg}}\label{THM1.5}
If MISP restricted to $H$-free tournaments can be solved in polynomial time, then $H$ is both 1-in-degenerate and 1-out-degenerate.
\end{customthm}
\proof{}  Let $\mathcal{C}$ be the class of graphs which can be obtained from replacing each edge of another graph by a path of 3 edges. By Lemma \ref{path3MVCP} and Theorem \ref{VC->MISP}, MISP on tournaments with backedge graph in $\mathcal{C}^+$ is NP-complete. We claim that if the backedge graph of a tournament $T$ is in $\mathcal{C}^+$, then $T$ is 1-in-degenerate and 1-out-degenerate. 

Let $G$ be a graph with $V(G) = \{v_1,...,v_n\}, E(G) = \{e_1,...,e_m\}$. Let $G'$ be the graph obtained by replacing each edge of $G$ with a path of 3 edges. We denote the two internal vertices corresponding to edge $e_i$ as $e_{ia}$ and $e_{ib}$. We now have $V(G') = \{v_1,...,v_n,e_{1a},e_{1b},...,e_{ma},e_{mb}\}$. We can now construct the graph ${G'}^+$ with vertex set $\{v_1,v_1',...,e_{mb},e_{mb}'\}$ and edge set $E(G')$. Let $T$ be the tournament whose backedge graph is ${G'}^+$ under enumeration $\sigma =(v_1,v_1',...,v_n,v_n',e_{1a},e_{1a}',e_{1b},...,e_{mb},e_{mb}').$
\begin{center}
\begin{tikzpicture}[thick,scale=0.4]
\node[draw,circle,inner sep=3pt,minimum size=3.5pt,label={$v_i$}] (1) at (-15,0){}; 
\node[draw,circle,inner sep=3pt,minimum size=3.5pt,label={$v_i'$}] (2) at (-12,0){}; 
\node[circle,inner sep=3pt,minimum size=3.5pt,label={$......$}] (x) at (-9,-1){};
\node[draw,circle,inner sep=3pt,minimum size=3.5pt,label={$v_j'$}] (3) at (-6,0){}; 
\node[draw,circle,inner sep=3pt,minimum size=3.5pt,label={$v_j'$}] (4) at (-3,0){}; 
\node[circle,inner sep=3pt,minimum size=3.5pt,label={$......$}] (x) at (0,-1){}; 
\node[draw,circle,inner sep=3pt,minimum size=3.5pt,label={$e_{ma}$}] (5) at (3,0){}; 
\node[draw,circle,inner sep=3pt,minimum size=3.5pt,label={$e_{ma}'$}] (6) at (6,0){}; 
\node[draw,circle,inner sep=3pt,minimum size=3.5pt,label={$e_{mb}$}] (7) at (9,0){}; 
\node[draw,circle,inner sep=3pt,minimum size=3.5pt,label={$e_{mb}'$}] (8) at (12,0){}; 
\draw[semithick,bend left = 35] (1) edge (5);
\draw[semithick,bend left = 25] (7) edge (5);
\draw[semithick,bend left = 40] (3) edge (7);
\end{tikzpicture}\\
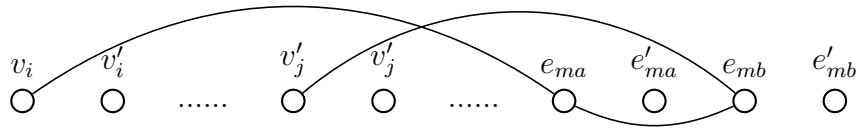
Figure 5: Example of a path of length 3 $(v_i,e_{ma},e_{mb},v_j)$ in $G^+$
\end{center}
Consider a subtournament of $T$ with vertex set $S$, and suppose for contradiction that $T[S]$ has no vertex with out-degree 1. Let $v$ be the last vertex of $S$ under enumeration $\sigma$. If $v \in \{v_1,v_1',...,v_n,v_n'\}$, then $v$ has outdegree 0. If $v \in \{e_{1a}',e_{1b}',...,e_{ma}',e_{mb}'\}$, then $v$ is not incident with any backedges, so $v$ still has outdegree $0$. If $v \in \{e_{1a},...,e_{ma}\}$ then $v$ is incident with $1$ backedge with vertices before $v$, so $v$ has outdegree at most $1$. It follows that $v \in \{e_{1b},...,e_{mb}\}$; let $v = e_{xb}$ for some $x \in [m]$. If $v$ has out-degree at least 2, then since $v$ is incident with at least 2 backedges in $T[S]$, it follows that $e_{ma} \in S$. If $e'_{ma} \in S$, $e_{ma}'$ has outdegree at most 1; if $e_{ma}' \notin S$, then $e_{ma} \in S$ has outdegree at most 1. Thus $T$ is 1-out-degenerate.

Now assume $H$ is not 1-out-degenerate, then it has a subtournament $H'$ in which all vertices have outdegree at least 2. By the above result, $T$ is $H'$-free, so $T$ is $H$-free. Moreover, by Lemma \ref{path3MVCP} and the reduction in Theorem \ref{npc}, we can use a solution for MISP in $T$ to solve MVCP in $G$. Thus, MISP for $H$-free tournaments is NP-complete. By reversing all edges of a tournament $G$, we get the complement $\overline{G}$. It is obvious that MISP for $G$ is equivalent to MISP for $\overline{G}$, so the above results also holds if we replace  `outdegree' by `indegree'. Hence if MISP for $H$-free tournaments can be solved in polynomial time, then $H$ is 1-in-degenerate and 1-out-degenerate. This concludes the proof. $\qed$

\begin{customcor}{\ref{t5f}}
MISP for $T_5$-free tournaments is NP-complete.
\end{customcor}
\subsection{Excluding $k$-snakes}
Given Theorem \ref{THM1.5}, it is natural to ask whether the converse holds. Namely, is it true that for every $H$ that is 1-in-degenerate and 1-out-degenerate, MISP for $H$-free tournaments can be solved in polynomial time? Unfortunately, the answer is NO. Here we give a tournament $H$ which is 1-in-degenerate and 1-out-degenerate, but we can prove MISP for $H$-free tournament is NP-complete.

\begin{sdef}
    A $k$-snake is the tournament with vertices $v_1,...,v_k$ such that the backedge graph $B_\sigma$ under the enumeration $\sigma=(v_1,...,v_k)$ has edge set $E(B_\sigma) = \{v_{i}v_{i+1}: i \in [k-1]\}$.
\end{sdef}

\begin{stheorem}
    A $k$-snake is both 1-in-degenerate and 1-out-degenerate.
\end{stheorem}
\proof{} Let $T$ be a $k$-snake and $S \subseteq V(T)$. Then, the backedge graph of $T[S]$ under the ordering $(v_1,...v_k)$ is a collection of paths. The first vertex and the last vertex in the backedge graph is incident with at most 1 backedges. Hence, the first vertex has in-degree at most 1 and the last vertex has out-degree at most 1. $\qed$

\begin{stheorem}
    MISP for 7-snake free tournaments is NP-complete.
\end{stheorem}
\proof{} (Sketch)
We use the same construction as in proof of Theorem \ref{THM1.5}, except that we need to change the enumeration order of $V(G'^+)$. For any edge $v_iv_j \in E(G)$ (assuming $i < j$), we denote the edge as $e_{ij}$. We order the edges such that $e_{ij,a}$ comes after $e_{kl,a}$ if and only if $i <k$ or $i=k,j< l$. Then we proceed as in the proof of Theorem \ref{THM1.5}.

Let $V = \{v_1,v_1',...,v_n,v_n'\}, E = V(G^+) \setminus V $. Notice that a 4-snake does not exist in the graph $G^+[E]$. Suppose there exists a 7-snake in $G^+$ with vertices $\{x_1,...,x_7\}$ such that under the ordering $(x_1,...,x_7)$, $x_{i}x_{i+1}$ is an edge for $i \in [6]$. According to whether $x_1,x_2,x_3$ are in $V$ or not, one can split into case studies. One can check in all cases, there is no candidate available in $G^+$ for the last vertex $x_7$, which shows that the graph $G^+$ is 7-snake free. Hence, MISP for 7-snake free tournaments is NP-complete. $\qed$
\section{Acknowledgement}
We are thankful to Sepehr Hajebi for many helpful discussions, and to Seokbeom Kim for allowing us to include the $(X, Y, Z)$-partition algorithm. 
\newpage
\printbibliography
\end{document}